\documentclass[10pt]{amsart}
\usepackage{latexsym,amssymb,graphicx}

\newtheorem{thm}{Theorem}[section]
\newtheorem{lem}[thm]{Lemma}

\newtheorem{cor}[thm]{Corollary}
\newtheorem{prop}[thm]{Proposition}
\theoremstyle{definition}

\newcommand{\comment}[1]{\marginpar{\sffamily{\tiny #1
\par}\normalfont}}
\renewcommand{\comment}[1]{}

\newcommand{\Z}{\ensuremath{{\mathbb{Z}}}}

\begin{document}

\title[]{On automorphism groups of free products of finite
groups, I: Proper Actions}

\author[]{Yuqing Chen}
      \address{Department of Mathematics\\
               RMIT University\\
               GPO Box 2476V\\
               MELBOURNE 3001, Australia}
      \email{rmayc@gauss.ma.rmit.edu.au}

\author[]{Henry H. Glover}
      \address{Department of Mathematics\\
               The Ohio State University\\
               Columbus, OH 43210}
      \email{glover@math.ohio-state.edu}

\author[]{Craig A. Jensen}
      \address{Department of Mathematics\\
               University of New Orleans\\
               New Orleans, LA 70148}
      \email{jensen@math.uno.edu}

\subjclass{20E36 (20J05)}
\keywords{automorphism groups, free products, proper actions}
\date{\today}

\begin{abstract}
If $G$ is a free product of finite groups, let $\Sigma Aut_1(G)$
denote all (necessarily symmetric) automorphisms of $G$ that do 
not permute factors in the
free product.
We show that a 
McCullough-Miller 
[D. McCullough and A. Miller,
{\em Symmetric Automorphisms of Free Products},
Mem. Amer. Math. Soc. 122 (1996), no. 582]
and  Guti\'{e}rrez-Krsti\'{c}
[M. Guti\'{e}rrez and S. Krsti\'{c},
{\em Normal forms for the group of basis-conjugating 
automorphisms of a free group}, 
International Journal of Algebra and Computation 8 (1998) 631-669]
derived
(also see Bogley-Krsti\'{c} 
[W. Bogley and S. Krsti\'{c},
{\em String groups and other subgroups of $Aut(F_n)$},
preprint]
space of pointed trees 
is an $\underline{E} \Sigma Aut_1(G)$-space for these groups.
\end{abstract}

\maketitle

\section{Introduction}\label{intro}

We remind the reader (see \cite{km,m}) that if $G$ is a 
discrete group, then the contractible $G$-space 
$\underline{E}G$ is characterized
(up to $G$-equivariant homotopy) by the property that
if $H$ is any subgroup of $G$ then the fixed point subcomplex
$\underline{E}G^H$ is contractible if $H$ is finite and empty 
if $H$ is infinite. These spaces are basic tools in studying the
geometry of the group $G$.

Let $G$ be a free product of $n$ finite groups.  
We wish to construct an
$\underline{E} \Sigma Aut_1(G)$-space  based on
McCullough-Miller's \cite{mm} space of trees,
which uses rooted trees similar to those found in 
Guti\'{e}rrez-Krsti\'{c} \cite{gk}.  
Here $\Sigma Aut_1(G)$ is the kernel
of the projection $\Sigma Aut(G) \to \Sigma_n$.  
We will show:

\begin{thm} Let $G$ be a free product of finite
groups and let $\Sigma Aut_1(G) = Aut_1(G)$ denote all
automorphisms of $G$ that do not permute groups in the
free product.  The space $L(G) = \underline{E} \Sigma Aut_1(G).$
That is, $L(G)$ is a contractible space which
$\Sigma Aut_1(G)$ acts on with finite stabilizers and
finite quotient.  Moreover, if $F$ is a finite subgroup
of $\Sigma Aut_1(G)$, then the fixed point subcomplex
$L(G)^F$ is contractible.
\end{thm}

We conjecture that the space $L(G)$ is in fact an
$\underline{E} \Sigma Aut(G)$-space in addition to
being an $\underline{E} \Sigma Aut_1(G)$-space.
We pause to note a few other related papers.
In \cite{cz} Collins and Zieschang establish the peak
reduction methods that underly all of the contractibility
arguments here.  Gilbert \cite{g} further refines these
methods and gives a presentation for $\Sigma Aut(G)$.
In \cite{bm} Bridson and Miller show that every finite
subgroup of $\Sigma Aut_1(G)$ fixes a point of 
McCullough-Miller space $K_0(G)$.  In \cite{bk}
Bogley and Krsti\'{c} completely calculate the cohomology
of $\Sigma Aut(F_n)$.  Brady, 
McCammond, Meier, and Miller \cite{bmmm} use
McCullough-Miller space to show that $\Sigma Aut(F_n)$ is
a duality group.

The authors would like to thank Mike Davis for
finding an error in an earlier version of this 
paper.

The paper is structured as follows.  In section 2,
we define the space $L(G)$ and in section 3 we show it is
contractible using a standard norm.  
In section 4 we briefly examine stabilizers
of points in $L(G)$.  In section 5, we develop many new norms,
each of which can be used to show $L(G)$ is contractible.
In section 6, we classify fixed point subcomplexes
$L(G)^F$ where $F$ is a finite subgroup of $\Sigma Aut_1(G)$
and in section 7 we show that these subcomplexes are
contractible. 

\section{Preliminaries}\label{prelim}

If $G = G_1 \ast G_2 \ast \cdots \ast G_n$, set $\mathcal{J} = 
\{G_1, \ldots, 
G_n\}$ and $\mathcal{J}^{0} = \{*, G_1, \ldots, G_n\}$.
For each $i$, choose an $1 \not = \lambda_n \in G_i$.
  Let
$\mathcal{P}(\mathcal{J}^{0})$ be the Whitehead poset constructed 
in \cite{mm}.
Elements of $\mathcal{P}(\mathcal{J}^{0})$ correspond to labelled 
bipartite trees,
where the $n+1$ labels come from the set $\{*, G_1, \ldots, G_n\}$.  
Often, as in \cite{mm}, we will abuse notation and take the labels from 
the set $\{*, 1, 2, \ldots, n\}$. 
Given a labelled tree $T$ and a labelled vertex $k$ in the tree,
two other labels are equivalent if they are in the same
connected component of $T-\{k\}$.  This gives us a partition
$\underline{\underline{A}}(k)$ of $\{*, 1, 2, \ldots, n\}$.
The singleton set
$Q(\underline{\underline{A}}(k))= \{k\}$ is called the operative 
factor of the partition.  Denote by $\underline{\underline{A}}$
the collection of all of these partitions as $k$ ranges over
$\{*, 1, 2, \ldots, n\}$.  This yields an equivalent notion of
elements of $\mathcal{P}(\mathcal{J}^{0})$.
The poset structure in 
$\mathcal{P}(\mathcal{J}^{0})$ comes from an operation called folding
(when the elements are thought of as labelled trees) or
by setting $\underline{\underline{A}}(k) \leq \underline{\underline{B}}(k)$
if elements of $\underline{\underline{A}}(k)$ are
unions of elements of $\underline{\underline{B}}(k)$. See 
McCullough and Miller \cite{mm} for more details.

Form a deformation retract of $\mathcal{P}(\mathcal{J}^{0})$ by folding 
all edges
coming in to $*$ on a labelled tree together, resulting in a labelled 
tree where $*$
is a valence 1 vertex.  Call the resulting poset $P(\mathcal{J})$.  
Observe that
elements of $P(\mathcal{J})$ correspond to pointed trees with labels 
in $\mathcal{J}$,
and that $P(\mathcal{J})$ is $(n-1)$-dimensional.  Denote elements of 
$P(\mathcal{J})$ as pairs $(\mathcal{J}, \underline{\underline{A}})$ as 
in \cite{mm}.

Now mimic the construction in section 2 of \cite{mm}.  That is, we 
must construct 
a space out of the posets $P(\mathcal{J})$.  Let $\mathcal{B}$ be the 
set of all bases 
of $G$.  Define a relation on 

$$\{(\mathcal{H},\underline{\underline{A}}) | \mathcal{H} \in \mathcal{B}, 
\underline{\underline{A}} \in P(\mathcal{H})\}$$

by relating $(\mathcal{H},\underline{\underline{A}})$ and
$(\mathcal{G},\underline{\underline{B}})$ whenever there is a product 
$\rho$ of symmetric Whitehead automorphisms carried by 
$(\mathcal{H},\underline{\underline{A}})$ so that $\rho \mathcal{H} = 
\mathcal{G}$
and $\rho \underline{\underline{A}} = \underline{\underline{A}}$. 
Denote the equivalence class by $[\mathcal{H},\underline{\underline{A}}]$
An 
automorphism
$(\mathcal{H},x)$ is {\em carried} by $(\mathcal{H},\underline{A})$ if 
$* \not = Q(\underline{A}) = Q(x)$, $x$ is constant on each
petal of $\underline{A}$, and $x$ is the identity on the petal containing $*$.
The set of all such equivalence classes forms a poset under the folding 
operation.  
Denote the poset $L(G)$.  

Recall that if all of the factors in $G$ are finite, 
the Kurosh subgroup theorem implies that $Aut(G)= \Sigma Aut(G)$.
Define an action of $\Sigma Aut(G)$ on $L(G)$ by having
$\phi \in \Sigma Aut(G)$ act via 
$\phi \cdot [\mathcal{H},\underline{\underline{A}}] = 
[\phi(\mathcal{H}),\phi(\underline{\underline{A}}))].$

Recall that the {\em symmetric Fouxe-Rabinovitch} subgroup 
$\Sigma FR(G)$ is the subgroup generated by all symmetric Whitehead 
automorphisms which do not conjugate their operative factor and that
$$\Sigma Aut_1(G) = \Sigma FR(G) \rtimes \Phi, 
\Sigma Aut(G) = \Sigma Aut_1(G) \rtimes \Omega$$
where $\Phi = \prod Aut(G_i)$ is the subgroup of factor automorphisms
and $\Omega$ (a product of symmetric groups) permutes the factors.
Further note that $\Phi$ and $\Omega$ are {\em not}
canonical.  Throughout this entire paper, 
We make the convention of choosing
them to be with respect to
the basis $\mathcal{H}_0 = \{G_1, \ldots, G_n\}$.

\section{Reductivity lemmas of McCullough and Miller}\label{reductivitymm}

In this section, we sketch how the work of 
McCullough and Miller in \cite{mm} implies that $L(G)$ is
contractible.  They show $K(G)$ is contractible
by defining a norm on nuclear vertices of $K(G)$
and inductively adding the stars of 
nuclear vertices using this norm while insuring that each new 
intersection in 
contractible.

We adopt an analogous approach.  First, we show that 
$L(G)$ is contractible using a norm which is directly analogous to
that of \cite{mm}.  In a later section, we will modify this norm
along the lines of Krsti\'{c} and Vogtmann in \cite{kv} and show that 
$L(G)$ can also be shown to be contractible with the modified norm. 

If $\mathcal{W}$ is a set of elements of $G$, we can define a norm on 
nuclear vertices
$[\mathcal{H}, \underline{\underline{0}}]$ of $L(G)$ by setting 
$\|\mathcal{H}\|_{\mathcal{W}}
=\sum_{w \in \mathcal{W}} |w|_{\mathcal{H}}$, where $|w|_{\mathcal{H}}$ is 
the
(non-cyclic) word length of $w$ in the basis $\mathcal{H}$.  

To avoid re-doing work that McCullough and Miller have already done, we
adopt the following conventions.  
Let $G$ be a free product of $n$ finite 
groups, as already noted.  
Let $\bar{G} = G \ast \langle \lambda_{n+1} \rangle$,
where $\langle \lambda_{n+1} \rangle \cong \Z/2$.  There is an injective map
from $\Sigma Aut(G)$ to $\Sigma Out(\bar{G})$ by sending 
$\phi \in \Sigma Aut(G)$ to
$\bar{\phi} \in \Sigma Out(\bar{G})$ where $\bar{\phi} (\lambda_{n+1}) = 
\lambda_{n+1}$. 
Moreover,
if $v = [\mathcal{H}, \underline{\underline{A}}]$ is a vertex in
$L(G)$, we can construct a corresponding vertex $\bar{v}$
in $K(\bar{G})$ by adding $\lambda_{n+1}$ to $\mathcal{H}$ and relabelling
the vertex $*$ in the tree corresponding to $\underline{\underline{A}}$
as $\langle \lambda_{n+1} \rangle$ (or just $n+1$.)  Note that if 
$(\mathcal{H}, x)$ is carried by $[\mathcal{H}, \underline{\underline{A}}]$
then $(\bar{\mathcal{H}}, \bar{x})$ 
is carried by $(\bar{\mathcal{H}}, \bar{\underline{\underline{A}}})$.
Finally, if $\mathcal{W}$
is a set of words in $G$, we can construct $\bar{\mathcal{W}}$ by sending
$w \in \mathcal{W}$ to $\bar{w} = w \lambda_{n+1} \in \bar{\mathcal{W}}$ 
(cf. Proposition 2.18 in \cite{g},
which is the basic idea of what we are doing here in this adjustment.)
Then $|w|_v+1=|\bar{w}|_{\bar{v}}$ for any $w \in \mathcal{W}$ so that 
$red_{\mathcal{W}}(\alpha, v) = 
red_{\bar{\mathcal{W}}}(\bar{\alpha}, \bar{v})$ 
for any $\alpha \in \Sigma Aut(G)$.

\begin{thm} \label{thm:premmcontractible}
The space $L(G)$ is contractible.
\end{thm}

\begin{proof}  We sketch how the work in Chapters 3 and 4 of
McCullough and Miller still applies.
Extension of Lemma 3.1 is trivial.  For Lemma 3.2, observe that we can 
define refinements and disjunctions for partitions of $\mathcal{J}^{0}$
as before, and that property 4 on page 27 of \cite{mm} implies that
if either $Q(\underline{A})=*$ or $Q(\underline{S})=*$, so that the 
relevant partition is trivial,
then the refinement or disjunction of $\underline{A}$ with 
$\underline{S}$ is just $\underline{A}$.  Hence analogs of Lemma 3.2 and
Lemma 3.3 follow in the context of $L(G)$.

Lemma 3.6 and 3.7 of \cite{mm} can be used to prove the analogous 
results in our new context.  (The reductive automorphism contructed in 
Lemma 3.7 still sends $\lambda_{n+1}$ to $\lambda_{n+1}$ because of 
the notions of
{\em constricted} peak reduction in \cite{g}.)

For Lemma 3.8 (which McCullough-Miller use to prove their Lemma 4.8 and at 
the end of their Lemma 4.9) we have symmetric Whitehead automorphisms
$\alpha, \sigma$ at $v$ in $L(G)$ and construct $\alpha_0, \sigma_0$ as in 
\cite{mm}.  The complication is that $\alpha_0, \sigma_0$ might
not be in $L(G)$ because the petal containing $*$ is not conjugated by 
the identity.  We resolve this by conjugating the whole 
automorphism, if necessary, so that the
petal containing $*$ does correspond to the identity.  
More specifically, construct the
corresponding $\bar{\alpha}, \bar{\sigma}, \bar{\alpha}_0, \bar{\sigma}_0$
in $K(\bar{G})$.  By McCullough-Miller's Lemma 3.8, we have
$$red_{\bar{\mathcal{W}}}(\bar{\alpha}_0, v) + 
red_{\bar{\mathcal{W}}}(\bar{\sigma}_0, v) \geq 
red_{\bar{\mathcal{W}}}(\bar{\alpha}_0, v) + 
red_{\bar{\mathcal{W}}}(\bar{\sigma}_0, v).$$
Also, 
$$red_{\bar{\mathcal{W}}}(\bar{\alpha}_0, v) + 
red_{\bar{\mathcal{W}}}(\bar{\sigma}_0, v) =
red_{\mathcal{W}}(\alpha_0, v) + red_{\mathcal{W}}(\sigma_0, v)$$
by our earlier observations.
Now conjugate (in $Aut(\bar{G})$) $\bar{\alpha}_0, \bar{\sigma}_0$ to 
obtain $\bar{\alpha}_0', \bar{\sigma}_0'$ which are the identity on
the petal containing $\lambda_{n+1}$.  Since these only differ by conjugation,
$$red_{\bar{\mathcal{W}}}(\bar{\alpha}_0, v) + 
red_{\bar{\mathcal{W}}}(\bar{\sigma}_0, v) =
red_{\bar{\mathcal{W}}}(\bar{\alpha}_0', v) + 
red_{\bar{\mathcal{W}}}(\bar{\sigma}_0', v)$$
and we can take the corresponding $\alpha_0', \sigma_0'$ in $L(G)$ so that
$$red_{\mathcal{W}}(\alpha_0', v) + red_{\mathcal{W}}(\sigma_0', v) \geq 
red_{\mathcal{W}}(\alpha_0, v) + red_{\mathcal{W}}(\sigma_0, v),$$
proving the analog of the lemma.

For Chapter 4, reason as follows.  Set 
$\mathcal{W}_0=\{\lambda_1, \ldots, \lambda_n\}$ so that the 
analog of Lemma 4.1 is
that there is only one nuclear vertex of minimal height $n$ in $L(G)$. 
Observe that every nuclear vertex in $K(\bar{G})$ corresponding to a 
basis of the form $\{\lambda_{n+1}^{i_1} \lambda_1 \lambda_{n+1}^{-i_1}, \ldots, 
\lambda_{n+1}^{i_n} \lambda_n \lambda_{n+1}^{-i_n}, \lambda_{n+1}\}$ is of minimal height $2n$
under the basis $\bar{\mathcal{W}}_0=\{\lambda_1 \lambda_{n+1}, \ldots, \lambda_n \lambda_{n+1}\}$.

Lemmas 4.2 and 4.3 are general lemmas about posets
from Quillen \cite{q} 
and hold without any modification.  
Theorem 4.4 and Proposition 4.5 are the central theorems, established in
section 4.2 by the lemmas 4.6, 4.7, 4.8, and 4.9.

Lemmas 4.6 and 4.7 can be proven using the same proof.  Lemma 4.8 can also
be proven using the same proof, even though it uses Lemma 3.8 which has been 
modified slightly.    The same holds for the crucial lemma, Lemma 4.9.  The
basic idea is that we could think of many of the calculations as taking place
in $K(\bar G)$, but just with Whitehead automorphisms whose domain
(cf. \cite{g} for notions of domain and constricted peak reduction) does 
not include $G_{n+1}$.  When we combine and modify these automorphisms,
we still obtain ones that are the identity on the last factor. 
\end{proof}

\section{Finite subgroups}\label{finite}

\begin{prop} \label{prop:finitestabilizers}
The stabilizer of a simplex in $L(G)$ under the
action of $\Sigma Aut(G)$ is finite.
\end{prop}

\begin{proof}
Let $[\mathcal{H}, \underline{\underline{A}}]$ be a vertex of $L(G)$.
If $i \not = j$, and both $(\mathcal{H}, x^i)$  and $(\mathcal{H}, x^j)$ are
symmetric Whitehead automorphisms carried by 
$[\mathcal{H}, \underline{\underline{A}}]$, then one of $x_j^i$ or $x_i^j$
must be the identity because at least one of
the petal of $\underline{\underline{A}}(i)$ containing $j$
or the petal of $\underline{\underline{A}}(j)$ containing $i$ also
contains $*$.
 Hence Lemma 7.4 of \cite{mm} applies to
give us that $(\mathcal{H}, x^i)$ and $(\mathcal{H}, x^j)$ commute.
Thus the stabilizer of 
$[\mathcal{H}, \underline{\underline{A}}]$ must be finite.
\end{proof}

\begin{prop} \label{prop:fixespoint}
Every finite subgroup of $\Sigma Aut_1(G)$ fixes a point
of $L(G)$.
\end{prop}

\begin{proof}
From Bridson and Miller in \cite{bm}, every finite subgroup of
$Aut(G)$ fixes a point
$v'=[\mathcal{H}', \underline{\underline{A'}}]$ of $K(G)$.
From Theorem 7.6 of \cite{mm}, any finite subgroup $F$ of 
$\Sigma Aut_1(G)$ that fixes $v'$
is conjugate by an inner automorphism $\mu$ 
to a subgroup whose elements are of the form
$$\prod (\mathcal{H}, x^i) \phi_i$$
where each $\phi_i \in Aut(G_i)$,the symmetric Whitehead automorphisms
$(\mathcal{H}, x^i)$ are carried by $v'$, and each $x_i^i=1$.
Moreover, there is a factor $k$ such that $x^i_k=1$ for all $i$ and there
is an unlabelled vertex $r$ of the tree $T'$ corresponding to 
$\underline{\underline{A'}}$ such that the petal containing $r$ is
always conjugated by the identity in the above Whitehead automorphisms.
Let $\mathcal{H} = \mu^{-1}(\mathcal{H}')$.  Form a tree $T$ by attaching
a free edge with terminal vertex $*$ to $T'$ at the vertex $r$, and let
$\underline{\underline{A}}$ be the vertex type determined by $T$.  
Then $F$ fixes the vertex $v=[\mathcal{H}, \underline{\underline{A}}]$
of $K(G)$ and every element of $F$ can be written in the form 
$$\prod (\mathcal{H}, x^i) \phi_i$$
where each $\phi_i \in Aut(G_i)$, the symmetric Whitehead automorphisms
$(\mathcal{H}, x^i)$ are carried by $v$, each $x_i^i=1$, and 
there is a factor $k$ such that $x^i_k=1$ for all $i$.
\end{proof}

\section{Better norm}\label{norm}

Well-order $G$ as $g_1, g_2, g_3, \ldots$ and order $\Z^G$ 
lexicographically.  For a nuclear vertex $v$ corresponding
to a basis $\mathcal{H}$, define a norm $\|v\| \in \Z^G$
by setting $\|v\|_i$ to be the (non-cyclic)
length $|g_i|_{\mathcal{H}}$
of $g_i$ in the basis $\mathcal{H}$.  This is analogous to 
the norm used by Krsti\'{c} and Vogtmann in \cite{kv} or Jensen
in \cite{j}.

\begin{prop} \label{prop:wellord}
The norm $\|\cdot\| \in \Z^G$ 
well orders the nuclear vertices of $L(G)$.
\end{prop}

\begin{proof}
Let $U$ be a nonempty subset of nuclear
vertices of $L(G)$ and proceed as in \cite{kv}.  That is,
inductively define $U_i$ and $d_i$ by setting $d_i$ to be
the minimal length $|g_i|_{\mathcal{H}}$ obtained by
all vertices $[\mathcal{H}, \underline{\underline{0}}] \in U_{i-1}$
and letting $U_i$ be all vertices of $U_{i-1}$ which obtain this
minimal length.  Recall that we chose specific $\lambda_i \in G_i$.
Let $N$ be such that $\lambda_i \in \{g_1, g_2, \ldots, g_N\}$ for
all $i$.  Let $\mathcal{W}_0=\{\lambda_1, \ldots, \lambda_n\}$.  We claim
that $U_N$ is finite.  Let $D=\sum_{g_{i_j}=\lambda_i} d_{i_k}$ so
that $\|v\|_{\mathcal{W}_0} \leq D$ for all $v \in U_N$.
Since each $G_i$ is finite, $L(G)$ is locally finite.
Hence the analog of the Existence Lemma 3.7 of \cite{mm} implies
the ball of radius $D$ (using the $\|\cdot\|_{\mathcal{W}_0}$
distance) around $\mathcal{H}_0$ in $L(G)$ is finite.  So
$U_N$ is finite.  Now choose $M \geq N$ large enough so that
$\{g_1, g_2, \ldots, g_M\}$ contains a representative from
each basis element of each basis corresponding to an element
of $U_N$.  Then $U_M$ contains exactly one element, 
the least element of $U$.
\end{proof}

Let $F$ be a finite subgroup of $\Sigma Aut_1(G)$.  Our goal in
the next few sections is to show that the fixed point subspace $L(G)^F$
is contractible.  A vertex $[\mathcal{H}, \underline{\underline{A}}]$
of $L(G)^F$ is {\em reduced} if no element of $L(G)^F$ lies below it in
the poset ordering.  We will show $L(G)^F$ is contractible by inductively
adding stars of reduced vertices and insuring that intersections are
always contractible.  The essential step will use the fact that
$L(G)$ can be shown the be contractible using the above norm, and we
will need the flexibility of being able to well-order $G$ in many
different ways.

\begin{thm} \label{thm:mmcontractible}
Given any well order of $G$, the norm $||v|| \in \Z^G$ defined
above on nuclear vertices of $L(G)$ is such that
$$st(v) \cap \left(\cup_{u < v} st(u)\right)$$
is contractible 
for any non-minimal nuclear vertex $v$, where $st(v)$ is the
star of $v$.  Hence $L(G)$ is contractible by induction.
\end{thm}

\begin{proof} We sketch how to apply \cite{mm} and 
Theorem \ref{thm:premmcontractible}.
To prove Lemmas 3.5, 3.6 of \cite{mm} with this new norm, simply apply
them, by letting $\mathcal{W}$ be a single word, in each coordinate
and applying the analogous lemmas
from Theorem \ref{thm:premmcontractible}.

To prove the Existence Lemma 3.7 of \cite{mm}, suppose $\mathcal{H}$ is
a given basis which does not have minimal norm.  Suppose that it does
have minimal norm on its first $m$ coordinates $\{g_1, g_2, \ldots, g_m\}$ 
but that the length of $g_{m+1}$ in $\mathcal{H}$ is not minimal.  Now
apply the Existence Lemma 3.7 of \cite{mm} with the set of words
defined to be $\{g_1, g_2, \ldots, g_m, g_{m+1}\}$ to get the desired
result.

For Lemma 3.8 (The Collins-Zieschang Lemma), note that the result is
proven in \cite{mm} by showing that the inequality holds coordinate-wise 
in our norm.

The arguments given in chaper 4 of \cite{mm} also carry through, except
that they are simplified somewhat because reductive edges now must be
strictly reductive. 
\end{proof}

\section{Fixed point subspaces.}\label{frfg}

If $[\mathcal{H},\underline{\underline{A}}]$ is a vertex type,
we think of $\underline{\underline{A}}$ as a
collection of partitions of $\{*, 1, 2, \ldots, n\}$ rather than a
collection of partitions of $\mathcal{H}$, where $wG_iw^{-1}$
in $\mathcal{H}$ is identified with $i$.
For each $k$, let $I_k(\underline{\underline{A}})$ be the 
set of labelled vertices that are a distance $2k$ away from $*$ in the
tree $T$ corresponding to $\underline{\underline{A}}$.  Let 
$I(\underline{\underline{A}}) = \cup_k I_k(\underline{\underline{A}})$
and define a poset order in $I=I(\underline{\underline{A}})$ 
by setting $r \leq s$ if $r$ occurs
on the minimal path in $T$ from $s$ to $*$.  For 
$i \in I_k=I_k(\underline{\underline{A}})$, define
$I(i)$ and $J(i)$ as follows.  Let $z_0=*, z_1, \ldots, z_k=i$ denote
the labelled vertices in the unique minimal path from $*$ to $i$ in $T$
and set $J(i)=(z_1, z_2, \ldots, z_k)$.  Let $a$ be the
unique unlabelled vertex between $z_{k-1}$ and $z_{k}$ and let
$I(i)=I(a)$ denote the set of all labels in $\{1,2, \ldots, n\}$ at a
distance $1$ from $a$ in $T$.  That is,
$I(i) = \{z_{k-1}\} \cup \{z \in I_k: z_{k-1} < z\}$.
(Exception: if $k=1$, let $I(i) = \{z \in I_k: z_{k-1} = * < z\}$.)  
Define $J_{<i}(\underline{\underline{A}})$ to be 
$\{z_1, z_2, \ldots, z_{k-1}\}$. Note that 
$J_{<i}(\underline{\underline{A}})$ is empty if 
$i \in I_1(\underline{\underline{A}})$.
Define words $w_i \in G$ inductively as follows.  For each $i \in I_1$, 
define $w_i \in G$ such that $H_i=w_i G_i w_i^{-1}$ and so that
$w_i$ has minimal length in the basis $\mathcal{H}$ (i.e., if
$H_i = w G_i w^{-1}$, then any word $wg_i$, $g_i \in G_i$ satisfies
$H_i = (wg_i) G_i (wg_i)^{-1}$ 
and we can choose one with minimal length.)  
For $i \in I_k$, 
$J(i)=(z_1, \ldots, z_k=i)$, define $w_i=w_{z_k}$ to be the word of minimal
length such that 
$$H_i = w_{z_1} w_{z_2} \ldots w_{z_k} G_i w_{z_k}^{-1} \ldots w_{z_1}^{-1}
w_{z_1}^{-1}.$$  For convenience, let 
$w(J(i)) = w_{z_1} w_{z_2} \ldots w_{z_k}$ so that 
$H_i = w(J(i)) G_i w(J(i))^{-1}$.

Let $a$ be an unlabelled vertex of the tree $T$ corresponding to 
$[\mathcal{H},\underline{\underline{A}}]$ which is at distance
$2k+1$ from $*$. 
If $k=0$, define the {\em stem} of $a$ to be $*$. 
If $k > 0$, define its {\em stem} to be the first labelled 
vertex on the unique shortest path from $a$ to $*$.
In either case, if $i$ is the stem of $a$ then define 
$\mathcal{H}(a) = \{w(J(i))^{-1} H_j w(J(i)) : j \in I(a)\}.$

For a given index $i \in I$, let $\pi_i: G \to G_i$ be the 
canonical projection.

\begin{lem} \label{lem:frfgtwisting}
Let
$[\mathcal{H}_0, \underline{\underline{A}}]$ be a vertex type and let
$\phi = \prod_j (\mathcal{H}_0, y^j) \psi_j$,
where each  $(\mathcal{H}_0, y^j)$ is symmetric
Whitehead automorphisms, $y_j^j=1$ for all $j$,
and each $\psi_j$ is a factor automorphism of $G_j$. 
Further suppose that $[\mathcal{H}_0,\underline{\underline{A}}]$ is
reduced in $L(G)^F$, where $F=\langle \phi \rangle$.  
Suppose
some other vertex type $[\mathcal{H},\underline{\underline{B}}]$
is also reduced in $L(G)^F$. 
Write $\phi = \prod_j (\mathcal{H}, x^j) \phi_j$
in this new basis, where $x_j^j=1$ for all $j$.
Write $H_i = w(J(i)) G_i w(J(i))^{-1}$ as above.
Then all of the following hold
\begin{enumerate}
\item For all $i,r$, $i \not = r$,
$$\pi_r(x_i^r) = 
\pi_r(\phi(w(J(i)))) y_i^r \pi_r(w(J(i))^{-1})$$
\item For all $i,j,r$, $j \not = r$,
if there exists a $g_r \in G_r$ such that
$y_i^r = \pi_r(\phi(g_r)) y_j^r \pi_r(g_r^{-1})$
and $y_j^r \not = 1$ then $y_i^r=y_j^r$.
\item $\underline{\underline{A}} = \underline{\underline{B}}$
(as partitions of $\{*, 1, 2, \ldots, n\}$.)
\end{enumerate}
\end{lem}

\begin{proof}
For a given index $i$, write 
$\lambda_i$ minimally in the basis $\mathcal{H}$ as
$$\lambda_i = a_r^{-1} \cdot \cdots \cdot a_1^{-1} \cdot (w \lambda_i w^{-1}) 
\cdot a_1 \cdot \cdots \cdot a_r,$$
where $w=w(J(i)(\underline{\underline{B}}))$.
Now $\phi= \prod_{j} (\mathcal{H}, x^j)$ sends $\lambda_i$ to
$$\phi(w^{-1})cw\psi_i(\lambda_i)w^{-1}c^{-1}\phi(w)$$
where 
$c \in \ast_{j \in J_{<i}(\underline{\underline{B}})} H_j$ 
comes from symmetric Whitehead
moves $(\mathcal{H}, x^j)$ conjugating $(w \lambda_i w^{-1})$. 
Similarly, $\phi= \prod_{j} (\mathcal{H}_0, y^j) \psi_j$
sends $\lambda_i$ to
$$d\psi_i(\lambda_i)d^{-1}$$
where 
$d \in \ast_{j \in J_{<i}(\underline{\underline{A}})} G_j$ 
comes from symmetric Whitehead
moves $(\mathcal{H}_0, y^j)$ conjugating $\lambda_i$.
So 
$$\phi(w^{-1})cw\psi_j(\lambda_i)w^{-1}c^{-1}\phi(w)=
d\psi_i(\lambda_i)d^{-1}$$
and there exists a $g_i \in G_i$ such that 
$\phi(w^{-1})cw = dg_i$. Thus
$$\pi_r(x_i^r) = 
\pi_r(\phi(w(J(i)))) y_i^r \pi_r(w(J(i))^{-1})$$
as desired.

By way of contradiction, suppose there exist
indices $i, j, r \in I(\underline{\underline{A}})$ 
and $g_r \in G_r$ such that
$y_i^r = \pi_r(\phi(g_r)) y_j^r \pi_r(g_r^{-1})$
and 
$y_i^r \not = y_j^r \not = 1$.
Let
$S_j = \{k \in I(\underline{\underline{A}}): y_j^r=y_k^r\}$.
Since the $S_j$ petal is not the identity petal,
we can conjugate it by $g_r$.
Let 
$\mathcal{H}'$ be the basis obtained from
$\mathcal{H}_0$ by conjugating all of the
$G_k$, $k \in S_j$, by $g_r$.
Then $[\mathcal{H}', \underline{\underline{A}}']$
=$[\mathcal{H}_0, \underline{\underline{A}}]$
and 
$\underline{\underline{A}} = \underline{\underline{A}}'$
as partitions of $\{*, 1, 2, \ldots, n\}$.
However the previous paragraph yields that
if $k \in S_j$ then $(y')_k^r = 
\pi_r(\phi(g_r)) y_j^r \pi_r(g_r^{-1})$.
In other words if $i \not = r$, we can combine the 
$S_i = \{k \in I(\underline{\underline{A}}): y_i^r=y_k^r, k \not = r\}$ 
and $S_j$
petals of
$\underline{\underline{A}}'(r)$.
(If $i=r$ then $1 = y_i^r = \pi_r(\phi(g_r)) y_j^r \pi_r(g_r^{-1})$.
In this case, we can combine the $S_j$ petal with the
petal containing $*$.)
This contradicts
the fact that 
$[\mathcal{H}_0, \underline{\underline{A}}]$
is reduced in $L(G)^F$.

Finally, we must show that 
$\underline{\underline{A}} = \underline{\underline{B}}$
as partitions of $\{*, 1, 2, \ldots, n\}$.
First, show this under the assumption that
for each $i$,
$I_{<i}(\underline{\underline{A}}) = 
I_{<i}(\underline{\underline{B}})$
(as partitions of $\{*, 1, 2, \ldots, n\}$.)
Suppose $x_j^r=x_i^r$ but $y_i^r \not = y_j^r \not = 1$.
Since $\pi_r(x_j^r)=\pi_r(x_i^r)$, (1) yields that
$y_i^r = \pi_r(\phi(w_i^{-1}w_j)) y_j^r \pi_r((w_i^{-1}w_j)^{-1})$
and we can apply (2) to get $y_i^r = y_j^r$.
As this is a contradiction, 
$y_j^r = y_i^r$ whenever $x_j^r=x_i^r$.
By symmetry, 
$x_j^r = x_i^r$ whenever $y_j^r=y_i^r$.
So $\underline{\underline{A}} = \underline{\underline{B}}$
in this case.

We are now ready to show 
$\underline{\underline{A}} = \underline{\underline{B}}$
in general. 
By way of contradiction, assume
$r \in J_{<j}(\underline{\underline{A}})$
but $r \not \in J_{<i}(\underline{\underline{B}})$.
Then $x_j^r=1$ and $y_j^r \not = 1$.  Let
$i \not = r$ be the index for which $y_i^r=1$.
(Abusing notation, we might have to take $i=*$.)
From (1), $y_i^r = 1 = \pi_r(x_j^r) = 
\pi_r(\phi(w(J(j)))) y_j^r \pi_r(w(J(j))^{-1})$.
By (2), $1=y_i^r=y_j^r$, which is a contradiction.
So $J_{<i}(\underline{\underline{A}}) \subset
J_{<i}(\underline{\underline{B}})$.
By symmetry, 
$J_{<i}(\underline{\underline{B}}) \subset
J_{<i}(\underline{\underline{A}})$
as well.  Now apply the previous case.
\end{proof}

\begin{lem} \label{lem:frfg}
Let
$[\mathcal{H}_0, \underline{\underline{A}}]$ be a 
vertex type
and let
$\phi = \prod_j (\mathcal{H}_0, y^j) \psi_j$,
where each  $(\mathcal{H}_0, y^j)$ is symmetric
Whitehead automorphisms, $y_j^j=1$ for all $j$,
and each $\psi_j$ is a factor automorphism of $G_j$. 
Further suppose that $[\mathcal{H}_0,\underline{\underline{A}}]$ is
reduced in $L(G)^F$, where $F=\langle \phi \rangle$.
Fix an index $k$ and let $d = \prod_{j \in J_{<k}} y_k^j$
(written so that if $j_1 < j_2$ in $J_{<k}$ 
then $y_k^{j_1}$ occurs before $y_k^{j_2}$ in the product.)
Then a word $w \in \ast G_j$ whose last letter is not in $G_k$ satisfies 
$\phi(w)=d w g_k d^{-1}$ for some $g_k \in G_k$ if and only 
if $w \in \ast_{j \in I(k)} G_{j,k}^{\circ}$
where the groups $G_{j,k}^{\circ}$ are defined by
$G_{j,k}^{\circ} = \{g \in G_j: y_k^j g (y_k^j)^{-1} = \psi_j(g)\}.$
\end{lem}

\begin{proof}
Let $i$ be the next labelled vertex
on a path from $k$ to $*$ in the tree $T$ corresponding to
$\underline{\underline{A}}$.
If $w \in \ast_{j \in I(k)} G_j$,
$\psi_j (w) = w$ for all $j \in I(k)-\{i\}$,
and $\psi_i (g_i) = y_k^i g_i (y_k^i)^{-1}$ 
for all $g_i \in G_i$ occurring in the normal form of $w$,
it is clear that $\phi(w)=d w d^{-1}$.
For the other direction, assume $\phi(w)=d w g_k d^{-1}$ and
suppose by way of contradiction that
$w \not \in \ast_{j \in I(k)} G_j.$
Let $S$ be the set of all indices $j$ for which an element of 
$G_j$ is a substring of $w$.  

\noindent {\bf Case 1: } {\em There is an index $r \in S-I(k)$ such that
$J_{<r} - J_{<k} \not = \emptyset$ or such that 
$J_{<r} - J_{<k} = \emptyset$ but $r \not \in J_{<k}$.} 
Choose $r$ satisfying the above condition to be maximal in the poset $I$.  
Now choose the first occurance $g_r$ of an element of $G_r$ in $w$ and write
$w=u_1g_ru_2$. 
Let $c = \prod_{j \in J_{<r}} y_r^j$.
Now $\phi(w)=du_1g_ru_2g_kd^{-1}=\phi(u_1)c\psi_r(g_j)c^{-1}\phi(u_2).$  
Because $r$ is
maximal in $I$ with the given condition, $\phi$ does not introduce any
more words from $G_r$ into $w$.  Hence we must have $du_1=\phi(u_1)c$
and $g_r = \psi_r(g_r)$.  If $s$ is the greatest index in $J_{<r}$,
then $y_r^s = \pi_s(\phi(u_1))^{-1} y_k^s \pi_s(u_1)$ and so 
by (2) of Lemma \ref{lem:frfgtwisting}, $y_r^s=y_k^s$.
This contradicts the fact that $r \not \in I(k)$.

\noindent {\bf Case 2: } {\em $S \subseteq I(k) \cup J_{<k}$.}
Let
$\psi = (\prod_{j \in I(k)-J_{<k}} \psi_j)
\cdot (\prod_{j \in J_{<k}} (\mathcal{H}_0, y^j) \psi_j)$,
ordered so that if $j_1 < j_2$ in $I$ then the automorphisms
with index $j_1$ are evaluated first (i.e., occur later in the
listing above.)
Note
that $\psi(w)=\phi(w)$.  Let $r$ be the least index in $S$
(least in the poset $I$) are let $g_r$ be the first
occurance of $G_r$ in $w$.  Write $w=u_1g_ru_2$.  Let $s$ be the
next labelled vertex on a path from $r$ to $k$ in $T$ and 
set $c = \prod_{j \in J_{<r}} y_s^j$.  After applying the first
$|J_{<r}|$ moves $(\mathcal{H}_0, y^j) \psi_j$ of $\psi$ to $w$, the result
is $cu_1g_ru_2c^{-1}$.  Moreover, the number of times an element of
$G_r$ occurs in the string $cu_1g_ru_2c^{-1}$ is the same as the 
number of times it occurs in $du_1g_ru_2g_kd^{-1}$.  After applying
the next move $(\mathcal{H}_0, y^r) \psi_r$ to $cu_1g_ru_2c^{-1}$,
we have $cy_s^ru_1(y_s^r)^{-1} \psi_r(g_r) y_s^r \ldots$.
Let $\psi'$ denote the last  $(|J_{<k}|-|J_{<s}|) + |I(k)|$
moves of $\psi$.  Then applying $\psi'$ gives us
$cy_s^r\psi'(u_1)(y_s^r)^{-1} \psi_r(g_r) y_s^r \ldots
=\phi(u_1) c \psi_r(g_r) y_s^r \ldots$.
Because applying $\psi'$ will not introduce any more elements of 
$G_r$, we have $\phi(u_1)c=du_1$.  This means that 
$1 = \pi_i(\phi(u_1)^{-1})y_k^i\pi_i(u_1)$ and thus $y_k^i=1$
by (2) of Lemma \ref{lem:frfgtwisting}.  This is a contradiction.

As we reached a contradiction in both cases,
$w \in \ast_{j \in I(k)} G_j.$  Since $w$ does not end in
an element of $G_k$, $g_k=1$ and $\phi(w)=dwd^{-1}$.
If $g_j \in G_j$ is a letter occuring in the normal form of $w$,
then $\psi_j(g_j)=g_j$ if $j \not = i$ and
$y_k^j g_j (y_k^j)^{-1} = \psi_j(g_j)\}$ if $j=i$.  Thus
$g_j \in G_{j,k}^{\circ}$ as desired.
\end{proof}

We apologize for the
confusing parentheses in $y_k^j g (y_k^j)^{-1}$ above,
which denotes conjugating $g$ by $y_k^j$.
Observe that if $y_k^j$ is in the center of $G_j$ 
(in particular, if $G_j$ is abelian) then
$G_{j,k}^{\circ} = \{g \in G_j: g = \psi_j(g)\}$.
Many of the arguments in this paper would be
simplified if we were only working with abelian 
factor groups.

\begin{prop} \label{prop:frfg}
Let
$v_0 = [\mathcal{H}_0, \underline{\underline{A}}]$ be a vertex type and let
$\phi = \prod_j (\mathcal{H}_0, y^j) \psi_j$,
where each  $(\mathcal{H}_0, y^j)$ is a symmetric
Whitehead automorphism, $y_j^j=1$ for all $j$,
and each $\psi_j$ is a factor automorphism of $G_j$. 
Further suppose that $[\mathcal{H}_0,\underline{\underline{A}}]$ is
reduced in $L(G)^F$, where $F=\langle \phi \rangle$.  
A necessary and sufficient condition
for any other vertex type $v$ to be
reduced in $L(G)^F$
is that it have a representative
$v=[\mathcal{H},\underline{\underline{B}}]$
where
$\underline{\underline{A}} = \underline{\underline{B}}$
(as partitions of $\{*, 1, 2, \ldots, n\}$) and that 
that 
when we write $H_i = w(J(i)) G_i w(J(i))^{-1}$ as above
we have 
$w_k \in \ast_{j \in I(k)} G_{j,k}^{\circ}$
where the groups $G_{j,k}^{\circ}$ are defined by
$G_{j,k}^{\circ} = \{g \in G_j: y_k^j g (y_k^j)^{-1} = \psi_j(g)\}.$
 Moreover, if
$v=[\mathcal{K},\underline{\underline{C}}]$ is any
other representative, then 
\begin{enumerate}
\item $\underline{\underline{C}} = \underline{\underline{B}}$
as partitions of $\{*, 1, 2, \ldots, n\}.$
\item We can get from $(\mathcal{K},\underline{\underline{C}})$
to $(\mathcal{H},\underline{\underline{B}})$ by a series of
moves conjugating petals $S$ of various
$\underline{\underline{C}}(i)$ by
$w(J(i)) \pi_i(w_k^{-1}) w(J(i))^{-1}$ where $k \in S$
(where the $w(J(i))$ and $w_k$ are taken with respect to 
the $(\mathcal{K},\underline{\underline{C}})$
representative of $v$.)
\end{enumerate}
\end{prop}

\begin{proof} 
For sufficiency, we note that it is a direct check to see that
$(\mathcal{H}_0, y^j) \psi_j \cdot v = v$ for all $j$ if $v$
is as described above.  So $\phi = \prod_j (\mathcal{H}_0, y^j) \psi_j$
fixes $v$ as well.

For necessity, suppose that $[\mathcal{H},\underline{\underline{B}}]$ 
is a reduced vertex in $L(G)^F$.  Then $\phi$ must fix this vertex
type, which means that
$[\mathcal{H},\underline{\underline{B}}]$
is stabilized by a product 
$$\prod_{j} (\mathcal{H}, x^j) \phi_j$$
which equals $\phi$, where $x_j^j=1$ for all $j$
and each $\phi_j$ is a factor automorphism of $H_j$.
By (3) of Lemma \ref{lem:frfgtwisting}, 
$\underline{\underline{A}} = \underline{\underline{B}}$
as partitions of $\{*, 1, 2, \ldots, n\}$.

We show that the $w_i$ have the desired properties by
inducting on the distance from $i$ to $*$ in $T$.  
If $i \in I_1$, write $w_i = \bar w_i \bar g_i$, where
$\bar w_i$ does not end in an element of $G_i$.
then $\phi(\bar w_i \lambda_i \bar w_i^{-1})
= \phi(\bar w_i) \psi_i(\lambda_i) \phi(\bar w_i)^{-1}
= \phi_i(\bar w_i \lambda_i \bar w_i^{-1}).$
Thus there exists a $g_i \in G_i$ such that
$\phi_i(\bar w_i) = (\bar w_i) g_i$.
By Lemma \ref{lem:frfg} $g_i=1$,
$\bar w_i \in \ast_{j \in I(i)} G_j$, and 
$\psi_j(\bar w_i)=\bar w_i$ for all $j$.
By way of contradiction, 
suppose $\bar g_i \not = 1$ but $\psi_i(\bar g_i) \not = \bar g_i$.
Write $\bar w_i=u_1 \ldots u_s$ in
normal form in the basis $\mathcal{H}$,
where each $u_j$ comes from an $H_{i_j}$ with 
$i_j \in I_1$.
Because $w_i = \bar w_i \bar g_i = (\bar w_i \bar g_i \bar w_i^{-1}) u_1 \ldots u_s$,
has length less than $s$, we have 
$u_1 = \bar w_i \bar g_i^{-1} \bar w_i^{-1}$.
Then
$$\bar g_i^{-1} w_i^{-1} u_2 \ldots u_s =
\psi_i(\bar g_i^{-1}) w_i^{-1} \phi(u_2 \ldots u_s)$$
Thus if $f = \prod_j \phi_j$ then
$$\phi(u_2 \ldots u_s) = f(u_2 \ldots u_s) = 
(w_i \psi_i(\bar g_i) \bar g_i^{-1} w_i^{-1}) u_2 \ldots u_s$$
which contradicts the fact that the length of 
$u_2 \ldots u_s$ in $\mathcal{H}$ is $s-1$.

For the inductive step, consider an index $k$ and let
$i$ be the next labelled vertex on a path from $k$ to $*$
in $T$.  Let $w=w(J(i))$.
As in the proof part (1) of
Lemma \ref{lem:frfgtwisting}, we have
$c \in \ast_{j \in J_{<i}(\underline{\underline{B}})} H_j$ 
coming from symmetric Whitehead
moves $(\mathcal{H}, x^j)$ conjugating $(w \lambda_i w^{-1})$
and 
$d \in \ast_{j \in J_{<i}(\underline{\underline{A}})} G_j$ 
coming from symmetric Whitehead
moves $(\mathcal{H}_0, y^j)$ conjugating $\lambda_i$
such that $cw = \phi(w)dg_i$, where $g_i=1$ by the induction
hypothesis.  Thus $\pi_i(w)=\pi_i(\phi(w)).$
As in the basis step of the induction,
write $w_k = \bar w_k \bar g_k$ where $\bar w_k$ does not
end in an element of $G_k$.
We have
$$c x_k^i w \bar w_k g_k = \phi(w \bar w_k) d y_k^i$$
for some $g_k \in G_k$.
It follows that 
$$\phi(\bar w_k) = (d w^{-1} x_k^i w) \bar w_k g_k (d y_k^i)^{-1}.$$
By (1) of Lemma \ref{lem:frfgtwisting} and the inductive
hypothesis, 
$$\pi_i(x_k^i) = \pi_i(\phi(w \bar w_k)) y_k^i \pi_i(w \bar w_k)^{-1}
=\pi_i(w) \pi_i(\phi(\bar w_k)) y_k^i \pi_i(\bar w_k)^{-1} \pi_i(w)^{-1}.$$
If $\pi_i(\phi(\bar w_k)) y_k^i \pi_i(\bar w_k)^{-1} \not = y_k^i$, change
$\mathcal{H}$ by conjugating the entire petal $S$ of
$\underline{\underline{B}}$ containing $k$ by
$w \pi_i(\bar w_k^{-1}) w^{-1}$.  Then in the new 
$\mathcal{H}$,
$\pi_i(x_j^i) 
=\pi_i(w) y_k^i \pi_i(w)^{-1}$
for all $j$ in $S$.  
Hence $w^{-1} x_k^i w = y_k^i$ and so
$$\phi(\bar w_k) = (d y_k^i) \bar w_k g_k (d y_k^i)^{-1}.$$
Applying Lemma \ref{lem:frfg},
$g_k=1$ and
$\bar w_k \in \ast_{j \in I(k)} G_{j,k}^{\circ}$.
For reasons similar to the base case of the induction,
$\psi_k(\bar g_k) = \bar g_k$ as well.
\end{proof}

\begin{cor} \label{cor:frfg}
Let
$[\mathcal{H}_0, \underline{\underline{A}}]$ be a vertex type and let
$\phi = \prod_j (\mathcal{H}_0, y^j) \psi_j$,
where each  $(\mathcal{H}_0, y^j)$ is a symmetric
Whitehead automorphism, $y_j^j=1$ for all $j$,
and each $\psi_j$ is a factor automorphism of $G_j$. 
Further suppose that $[\mathcal{H}_0,\underline{\underline{A}}]$ is
reduced in $L(G)^{F_1}$, where $F_1=\langle \phi \rangle$.  
Let $F_2$ be the subgroup generated by all of 
the $(\mathcal{H}_0, y^j) \psi_j$, so
that $F_1 \subseteq F_2$.
Then some other  vertex type $[\mathcal{H},\underline{\underline{B}}]$
is reduced in $L(G)^{F_1}$ if and only if
it is reduced in $L(G)^{F_2}$.
\end{cor}

\begin{proof}
Since $F_1 \subseteq F_2$, $L(G)^{F_2} \subset L(G)^{F_1}$.
However, if $v$ is a nuclear vertex of $L(G)^{F_1}$,
then from Proposition \ref{prop:frfg} every element of
$F_2$ fixes it (i.e,  
$G_{j,k}^{\circ} = \{g \in G_j: y_k^j g (y_k^j)^{-1} = \psi_j(g)\}$
only depends on $(\mathcal{H}_0, y^j) \psi_j$.)
\end{proof}

Let
$v_0 = [\mathcal{H}_0, \underline{\underline{A}}]$ be a vertex type,
$F$ be a finite subgroup of $\Sigma Aut_1(G)$ that fixes $v_0$,
and suppose $v_0$ is reduced in $L(G)^F$.
Let $\phi = \prod_j (\mathcal{H}_0, y^j) \psi_j \in F$,
where each  $(\mathcal{H}_0, y^j)$ is symmetric
Whitehead automorphisms, $y_j^j=1$ for all $j$,
and each $\psi_j$ is a factor automorphism of $G_j$. 
Define $\pi_j(\phi) = (\mathcal{H}_0, y^j) \psi_j.$
Define
the groups $G_{j,k}^{\circ}$ by
$$G_{j,k}^{\circ} = \cap_{\phi \in F}
\{g \in G_j: y_k^j g (y_k^j)^{-1} = \psi_j(g) \hbox{ where } \pi_j(\phi) = 
(\mathcal{H}_0, y^j) \psi_j\}.$$
A representative $(\mathcal{H},\underline{\underline{B}})$
of some other 
vertex $v=[\mathcal{H},\underline{\underline{B}}]$ is
{\em $F$-standard} if all of the following hold
\begin{itemize}
\item $\underline{\underline{A}} = \underline{\underline{B}}$
as partitions of $\{*, 1, 2, \ldots, n\}$.
\item When we write $H_i = w(J(i)) G_i w(J(i))^{-1}$ 
then $w_i \in \ast_{j \in I(i)} G_{j,k}^{\circ}$.
\end{itemize}

\begin{thm} \label{thm:frfg}
Let
$[\mathcal{H}_0, \underline{\underline{A}}]$ be a vertex type 
and suppose that $F \subseteq \Sigma Aut_1(G)$ fixes
$[\mathcal{H}_0, \underline{\underline{A}}]$.
Further suppose that $[\mathcal{H}_0,\underline{\underline{A}}]$ is
reduced in $L(G)^F$.  
A necessary and sufficient condition
for any other vertex type $v$
to be reduced is that it have an
$F$-standard representative.
\end{thm}

\begin{proof}
Let $F_{+}$ be the subgroup generated by 
$\{\pi_j(\phi): j \in I, \phi \in F\}$
so that 
$F \subseteq F_{+}$. 
Now $L(G)^{F_{+}} = L(G)^{F}$ by Corollary \ref{frfg}.
From Proposition \ref{prop:frfg}, each reduced vertex 
$v=[\mathcal{H},\underline{\underline{B}}]$ in
$L(G)^{F_{+}}$ must have 
$\underline{\underline{A}} = \underline{\underline{B}}$
because the structure of $\underline{\underline{B}}$
depends only on the symmetric Whitehead moves
$(\mathcal{H}_0, y^j)$ occuring in $\phi \in F$.  In particular,
$\underline{\underline{B}}(j)$ is the wedge (see \cite{mm} page
29) of all of the full carriers of the
$(\mathcal{H}_0, y^j)$ occuring in $\phi \in F$.

To show that each $w_i \in \ast_{j \in I(i)} G_{j,k}^{\circ}$,
note that (2) of Proposition \ref{prop:frfg} means that if
the letters from $G_j$ in a particular $w_k$ are
already in 
$$\cap_{\phi \in F_1}
\{g \in G_j: y_k^j g (y_k^j)^{-1} = \psi_j(g) \hbox{ where } \pi_j(\phi) = 
(\mathcal{H}_0, y^j) \psi_j\},$$
then if we take a $\xi \not \in F$ and conjugate petals again to get the 
letters in
$$\{g \in G_j: y_k^j g (y_k^j)^{-1} = \psi_j(g) \hbox{ where } \pi_j(\xi) = 
(\mathcal{H}_0, y^j) \psi_j\},$$
then they are still in the previous group and hence in the intersection
of the two groups.
\end{proof}

\section{Contractibility of fixed point subspaces.}\label{contractible}

For this entire section, 
let $F$ be a finite subgroup of $\Sigma Aut_1(G)$ that fixes
$[\mathcal{H}_0, \underline{\underline{A}}]$ and suppose that
$[\mathcal{H}_0, \underline{\underline{A}}]$ is reduced in $L(G)^F$.
If $[\mathcal{H}, \underline{\underline{B}}]$ is any other reduced
vertex and $\mathcal{H} = \{H_1, \ldots, H_n\}$ then
for all $j,k$ define $H_{j,k}^{\circ} = w(J(i)) G_{j,k}^{\circ} w(J(i))^{-1}.$
If $a$ is the next unlabelled vertex on a path from $k$ to $*$,
set $G_{j,a}^{\circ}=G_{j,k}^{\circ}$
and $H_{j,a}^{\circ}=H_{j,k}^{\circ}$.  
In addition, define $G_{j,k}^{\circ\circ}$ to be 
$G_{j,k}^{\circ}$ if it is nontrivial and $\Z/2$ otherwise.
Define $H_{j,k}^{\circ\circ}$, etc., analogously. For each $j,a$
choose $1 \not = \lambda_{j,a} \in G_{j,a}^{\circ\circ}.$
Let $G_a = \ast_{j \in I(a)} G_{j,a}^{\circ\circ}.$  

Note that if $[\mathcal{H}, \underline{\underline{A}}] \in L(G)^F$
with $(\mathcal{H}, \underline{\underline{A}})$ $F$-standard
and $j \in I_1$, then 
$\ast_{k \in I_1} H_k = \ast_{k \in I_1} G_k$.
This follows by letting $N$ be the normal closure of 
$\ast_{k \not \in I_1} G_k$ and considering $G/N \cong \ast_{k \in I_1} G_k$.
Observe that $\ast_{k \not \in I_1} H_k \subseteq N$.
If $j \in I_1$, then $H_j \subseteq \ast_{k \in I_1} G_k$ because
$H_j = w_j G_j w_j^{-1}$ with $w_j \in \ast_{k \in I_1} G_k.$  Therefore,
$\ast_{k \in I_1} H_k \subseteq \ast_{k \in I_1} G_k$.  Furthermore,
since $\mathcal{H}$ is a basis of $G$, if $j \in I_1$ and $g_k \in G_k$,
then we can write $g_k = v_1 v_2 \ldots v_s$ in the basis $\mathcal{H}$.
Taking the quotient by $N$, this yields a way of writing $g_k$ in
$\ast_{k \not \in I_1} H_k$.  It follows that 
$\ast_{k \in I_1} H_k = \ast_{k \in I_1} G_k$, as desired.
More generally, one can verify that 
$\ast_{k \in I(a)} H_k = \ast_{k \in I(a)} w(J(i)) G_k w(J(i))^{-1}$.
A direct induction argument now yields that 
$w(J(i)) \in \ast_{k \in I(a) \cup J_{<i}} G_k$.

With the same setup and hypothesis of Theorem \ref{thm:frfg}
and where $[\mathcal{H}, \underline{\underline{A}}] \in L(G)^F$
with $(\mathcal{H}, \underline{\underline{A}})$ $F$-standard, we have

\begin{lem} \label{lem:frfgwordlength}
Let $a$ be an unlabelled vertex of $T$ with stem $i$.
Then for each $h \in \ast_{k \in I(a)} G_k$,
$$|h|_{\mathcal{H}} = 2|w(J(i))|_{\mathcal{H}} + 
|h|_{\mathcal{H}(a)}.$$
\end{lem}

\begin{proof}
Let $w = w(J(i))$.  We show the result by induction on the
distance $d$ from $a$ to $*$ in $T$.
Assume $d \geq 3$, as the basis step of $d=1$ is immediate.

If $m$ is the next labelled vertex on a path from $i$ to 
$*$ and $y$ is the unlabelled vertex between $i$ and $m$ then 
$w = w'w_i$, $w'=w(J(m))$.  By our inductive hypothesis,
and $|w_i|_{\mathcal{H}} = 2|w(J(m))|_{\mathcal{H}} + 
|w_i|_{\mathcal{H}(y)}.$  In other words, if 
$w' = v_{t+1} \ldots v_{s}$ in $\mathcal{H}$ 
and $w_i = (w'^{-1} v_{1} w') (w'^{-1} v_{2} w') 
\ldots (w'^{-1} v_{t} w')$ is a minimal way of writing
$w_i$ in $\mathcal{H}(y)$, then
$w_i = 
v_{s}^{-1} \ldots v_{t+1}^{-1} 
(v_{1} v_{2} \ldots v_t) v_{t+1} \ldots v_s$
is a minimal way of writing the length $t+s$ word $w_i$
and has no cancellations in $\mathcal{H}$.
Thus $w = w' w_i = (v_1 v_2 \ldots v_t) v_{t+1} \ldots v_s$
is a minimal way of writing the length $s$ word $w$.

If $h = u_1 u_2 \ldots u_r$ is a minimal way of writing 
$h$ in the basis $\mathcal{H}(a)$ then
$$h = v_s^{-1} \ldots v_1^{-1} (w u_1 w^{-1}) \ldots (w u_r w^{-1})
v_1 \ldots v_s.$$
No cancellation occurs among the $w u_j w^{-1}$ by themselves or the
$v_j$ by themselves.  We must verify that 
no cancellation occurs at the stages 
$v_1^{-1} (w u_1 w^{-1})$ or $(w u_r w^{-1}) v_1$ because
the $v_i$ are not in $H_l$ for $l \in I(a)$.  This follows
because 
if $v_1=w u_1 w^{-1}$, then $v_1 \in H_i$, $u_1 \in G_i$.
But recall that we chose $w_i$ to have minimal length
among all $w_i g$, $g \in G_i$, and 
$$v_{s}^{-1} \ldots v_{t+1}^{-1} 
(v_{2} \ldots v_t) v_{t+1} \ldots v_s= w_i u_1^{-1}$$ 
has smaller length.
So this does not occur and similarly no
other cancellations occur. 
\end{proof}

Let $a$ be an unlabelled vertex of the tree $T$ corresponding
to some $v=[\mathcal{H}, \underline{\underline{A}}]$
with $(\mathcal{H}, \underline{\underline{A}})$ $F$-standard
and let $i$ be the stem of $a$.  

Let 
$$\mathcal{H}_a = 
\{w_j G_{j,a}^{\circ\circ} w_j^{-1}: j \in I(a)-\{i\}\} \cup \{G_{i,a}^{\circ\circ}\}
= \{w(J(i))^{-1} H_{j,a}^{\circ\circ} w(J(i)): j \in I(a)\}.$$
Thus $w_j G_{j,a}^{\circ\circ} w_j^{-1} \in \mathcal{H}_a$
iff $w_j G_j w_j^{-1} \in \mathcal{H}(a)$.  By Theorem
\ref{thm:frfg}, each $w_j \in G_a = \ast_{j \in I(a)} G_{j,a}^{\circ\circ}.$

Let $A$ be the set of unlabelled 
vertices in the tree $T$ 
corresponding to $\underline{\underline{A}}$.  Well order
$A$ so that if $a$ is on the unique shortest path from $b$ to 
$*$, then $a < b$.   
Choose a well order for each $G_a$ that
puts $\lambda_{i,a}$, where $i$ is the stem of $a$ first, 
then the other letters
$\lambda_{j,a}$, $j \in I(a)-\{i\}$, and finally all of the other words.
Well order $\cup_{a \in A} G_a$ so that:
(i) If $g < h$ in $G_a$ then $g < h$ in $\cup_{a \in A} G_a$;
and (ii) If $a < b$, $i$ is the stem of $a$, and $j$ is the
stem of $b$, then every element of $G_a-G_{i,a}^{\circ\circ}$ occurs before
every element of $G_b-G_{j,b}^{\circ\circ}$.

Order $\Z^{\cup_{a \in A} G_a}$
lexicographically and define a norm 
$$\|(\mathcal{H}, \underline{\underline{A}})\| \in 
\Z^{\cup_{a \in A} G_a}$$
on the $F$-standard pair
$(\mathcal{H}, \underline{\underline{A}})$ 
representing a nuclear vertex of $L(G)^F$
by setting the $g$th coordinate to be $|g|_{\mathcal{H}},$ 
the length of the word $g$ in the 
basis given by $\mathcal{H}$.  
As stated, this does not define a norm on nuclear
vertices of $L(G)^F$ because there is more than
one way to write the vertex type
$v=[\mathcal{H}, \underline{\underline{A}}]$
in an $F$-standard basis.
We solve this problem by definining 
$$\|v\| = 
\hbox{min}_{[\mathcal{K}, \underline{\underline{A}}]=v, 
(\mathcal{K}, \underline{\underline{A}}) \hbox{ $F$-standard}}
\|(\mathcal{H}, \underline{\underline{A}})\|.$$
Observe that given any
representative
$(\mathcal{H}, \underline{\underline{A}})$
there is an easy algorithm to construct the minimal
representative by proceeding inductively through
$A$.  If $a_0$ is the least element of $A$
(the vertex adjacent to $*$ in $T$), then we
cannot change the values in the range $G_{a_0}$
at all.  Supposing we have minimized all values
less than a particular $a \in A$, we let $i$
be the stem of $a$.  We now conjugate all of
$I(a)-\{i\}$ by a single element of 
$w(J(i)) G_{i,a}^{\circ} w(J(i))^{-1}$, if necessary, to reduce 
the norm restricted to $G_a$.  

\begin{prop} \label{prop:frfgwellord}
The norm 
$$\|v\| \in 
\Z^{\cup_{a \in A} G_a}$$
well orders the nuclear
vertices of $L(G)^F$.
\end{prop}

\begin{proof}
Let $U$ be a nonempty subset of nuclear
vertices of $L(G)$. 
Inductively define $U_g$ and $d_g$ by setting $d_g$ to be
the minimal length $|g|_{\mathcal{H}}$ obtained by
all vertices 
$[\mathcal{H}, \underline{\underline{A}}] \in \cap_{h < g} U_h$
and letting $U_g$ be all vertices of $\cap_{h < g} U_h$ which obtain this
minimal length.

We show by induction that if $a$ is an unlabelled vertex of $T$
then any
$\mathcal{H}, \mathcal{K} \in \cap_{g \in G_a} U_g$
satisfy $\mathcal{H}_{a} = \mathcal{K}_{a}$.
For the basis step, let $a_0$ denote the unlabelled
vertex adjacent to $*$ and note that 
$\|\cdot\| \in \Z^{G_{a_0}}$ well orders the
nuclear vertices of $L(G_{a_0})$ by 
Proposition \ref{prop:wellord}.

For the inductive step, consider an unlabelled
vertex $a \not = a_0$ of $T$ and suppose
$\mathcal{H}, \mathcal{K} \in \cap_{g \in G_a} U_g$.
By induction, for all $b < a$, 
$\mathcal{H}_{b} = \mathcal{K}_{b}.$  
In particular, if $i$ is the stem of $a$, then
$w(J_{\mathcal{H}}(i))=w(J_{\mathcal{K}}(i))$
and 
$|w(J_{\mathcal{H}}(i))|_{\mathcal{H}}=
|w(J_{\mathcal{K}}(i))|_{\mathcal{K}}$.  Now
use Lemma \ref{lem:frfgwordlength} and 
Proposition \ref{prop:wellord} applied to
$L(G_a)$ to get that 
$\mathcal{H}_{a} = \mathcal{K}_{a}$.

So $\mathcal{H}=\mathcal{K}$ and we are done. 
\end{proof}

Observe that $[\mathcal{H}_0, \underline{\underline{A}}]$
is the unique minimal vertex 
(with $(\mathcal{H}_0, \underline{\underline{A}})$ 
as its minimal $F$-standard representative) of $L(G)^F$ in the above norm.
A strictly reductive symmetric Whitehead move at a
non-minimal nuclear vertex 
$[\mathcal{H}, \underline{\underline{A}}]$
is a symmetric Whitehead automorphism
carried by some vertex type in the
ascending star of $[\mathcal{H}, \underline{\underline{A}}]$.

The well order in $\cup_{a \in A} G_a$ restricts to 
a well order on $G_a$ so that we have an induced norm in $\Z^{G_a}$ on
$L(G_a)$.  
Let $v=[\mathcal{H}, \underline{\underline{A}}]$, where
we assume for the remainder of the section that whenever we
write $v$ this way, we have chosen
$(\mathcal{H}, \underline{\underline{A}})$ to be 
a minimally-normed $F$-standard pair representing $v$ using the 
algorithm stated before Proposition \ref{prop:frfgwellord}.
If $\alpha=(\mathcal{H}_a, y^k)$ is a reductive 
Whitehead move at $v_a = [\mathcal{H}_a, \underline{\underline{0}}]$ in 
$L(G_a)$, then let $\alpha^a=(\mathcal{H}, x^k)$ be defined 
as follows:
If $j \leq i$ in $I$ or $j$ is not comparable with $i$ in $I$,
then define $x_j^k=1$.  On the other hand, if
$j \geq l$ for some $l \in I(a)-\{i\}$, then set
$x_j^k = w(J(i)) y_l^k w(J(i))^{-1}$. 
Suppose  $\alpha$ is carried by 
$[\mathcal{H}_a, \underline{\underline{B}}]$ and $T_a$ is the
tree for $\underline{\underline{B}}$. 
Define $T^a$ by 
first cutting out $a$ and all edges attached to $a$ from $T$, and then 
glueing $T_a$ in to the resulting hole, by attaching the vertex $j$
of $T_a$ to the vertex $j$ of $T-\{a\}$.  Let 
$[\mathcal{H}, \underline{\underline{B}}^a]$ be the vertex type 
corresponding to the tree $T^a$ and observe that
$\alpha^a$ is carried by $[\mathcal{H}, \underline{\underline{B}}^a]$.

\begin{lem} \label{lem:frfgsmalltolarge}
If $\alpha=(\mathcal{H}_a, y^k)$ is a reductive 
Whitehead move at $v_a = [\mathcal{H}_a, \underline{\underline{0}}]$ 
in 
$L(G_a)$ carried by $[\mathcal{H}_a, \underline{\underline{B}}]$
as described above, then $\alpha^a=(\mathcal{H}, x^k)$
is a reductive Whitehead move at 
$v=[\mathcal{H}, \underline{\underline{A}}]$ and is
carried by $[\mathcal{H}, \underline{\underline{B}}^a]$.
\end{lem}

\begin{proof}
Recall that $i$ is the stem of $a$.
Let $w(J(i))=w$.
Note that we can assume $y_i^k=1$ since
$|\lambda_{i,a}|_{v_a}=1$ is minimal.
By way of contradiction, suppose
the index $k$ is one for which  $G_{k,a}^{\circ}=\langle 1 \rangle.$
Recall that by Theorem \ref{thm:frfg},
$w_j \in \ast_{r \in I(a)} G_{r,a}^{\circ}$
for all  $j \in I(a)-\{i\}$.
So when we write each $w_j$ in the basis $\mathcal{H}_a$, we
need not use any letters from $w(J(i))^{-1} H_{k,a}^{\circ\circ} w(J(i))$.
In addition,
if $j \not = k$, when we write $\lambda_{j,a}$ in the basis 
$\mathcal{H}_a$, we need
not use any letter from $w(J(i))^{-1} H_{k,a}^{\circ\circ} w(J(i))$.  
If $j = k$, the normal form of $\lambda_{j,a}$ in the basis 
$\mathcal{H}_a$ uses exactly one
letter from $w(J(i))^{-1} H_{k,a}^{\circ\circ} w(J(i))$
and the rest from other elements of $\mathcal{H}_a$.

Let $g = u_{1, r_1} \ldots u_{t, r_t}$ 
be the normal form in the basis $v_a$ 
of some element of $\ast_{j \in I(a)} G_{j,a}^{\circ}$,
where $u_{1, r_j} \in G_{r_j,a}^{\circ}.$  Then the
normal form of $g$ in the basis $\alpha(v_a)$ is
the product
$$[(y_{r_1}^k)^{-1}]  
[(y_{r_1}^k u_{1,r_1} (y_{r_1}^k)^{-1}] 
[y_{r_1}^k (y_{r_2}^k)^{-1}] 
[y_{r_2}^k u_{2,r_2} (y_{r_2}^k)^{-1}]$$ 
$$\ldots
[y_{r_{t-1}}^k (y_{r_t}^k)^{-1}]
[y_{r_t}^k u_{t,r_t} (y_{r_t}^k)^{-1}] 
[y_{r_t}^k]$$
which has length greater than or equal to $t$.
Furthermore, the length is equal to $t$ only when
$y_{r_j}^k=1$ for every $j = 1, 2, \ldots t$.

Taking $g$ to be $\lambda_{j,a}$ for $k \not = j \in I(a)$,
we see that $\alpha$ cannot reduce any of the
lengths $|\lambda_{j,a}|$  Taking $g$ to be $w_k$, 
we also see that $\alpha$ cannot reduce $|\lambda_{k,a}|$.
But $\alpha$ is reductive by hypothesis and the
first coordinates of $G_a$ are the $|\lambda_{j,a}|$
for $j \in I(a)$.  By the above paragraph,
each $y_j^k=1$ and $\alpha$ is the identity map.
This is a contradiction.  So 
$G_{k,a}^{\circ}$ is nontrivial.

Since
$\alpha^a$ was defined by letting $x_j^k=1$ if
$j \leq i$ in $I$ or $j$ is not comparable with $i$ in $I$,
we know that $\alpha^a$ does not change any 
coordinate $|g|_{\mathcal{H}}$ with $g \in G_b$ for some
$a \not = b \in X$ where (i) $b$ is on the path from $i$ the $*$ in $T$
or (ii) where $a$ is not on the path from $b$ to $*$. 
Let $h \in G_a$ give the first coordinate where
$\alpha$ is reductive.
By Lemma \ref{lem:frfgwordlength}, $\alpha^a$ is not reductive
on any coordinate of $\cup_{b \in X} G_b$ up to $h$, and it
is as reductive as $\alpha$ is on the $h$ coordinate.

For a particular $l \in I(a)-\{i\}$,  
$x_j^k$ is constant on the branch of $T$ given by taking 
$j \geq l$.  Hence $\alpha^a$ is carried by 
$[\mathcal{H}, \underline{\underline{B}}^a]$ in the ascending star of 
$v=[\mathcal{H}, \underline{\underline{A}}]$. 
\end{proof}

\noindent {\bf Note: } Observe that there are some cases
where a non-reductive move $\alpha=(\mathcal{H}_a, y^k)$ 
at $v_a$ still induces
a well defined (but non-reductive) move $\alpha^a$ at $v$.
In particular, we must have $G_{k,a}^{\circ}$ nontrivial
and $y_i^k=1$.

Next we investigate how a reductive move at $v$ 
defines moves at various $v_a$s.
Let $\alpha = (\mathcal{H}, y^k)$ be a reductive move at 
$v=[\mathcal{H}, \underline{\underline{A}}]$ 
carried by $[\mathcal{H}, \underline{\underline{B}}]$ and 
let $a$ be
an unlabelled vertex of $T$ which is adjacent to $k$.  
Suppose that $Y$ is the tree for $\underline{\underline{B}}$
We 
can define the move $\alpha_a= (\mathcal{H}_a, x^k)$ at
$v_a=[\mathcal{H}_a, \underline{\underline{0}}]$
carried by $[\mathcal{H}, \underline{\underline{B_a}}]$
as follows:

\noindent {\bf Case 1: } {\em $a$ is the next vertex on a path from
$k$ to $*$ in $T$.}  Let $i$ be the stem of $a$.  Set
$x^k = w(J(i))^{-1} y^k w(J(i)))$.  The tree $T_a$ for
$\underline{\underline{B_a}}$ is given by looking at the 
subtree of $Y$ spanned by vertices in $I(a)$.

\noindent {\bf Case 2: } {\em $k$ is the stem of $a$.}
Set $x^k = w(J(i))^{-1} y^k w(J(i)))$.
As in the previous case, the tree $T_a$ for
$\underline{\underline{B_a}}$ is given by looking at the 
subtree of $Y$ spanned by vertices in $I(a)$.

\begin{lem} \label{lem:frfglargetosmall}
If $\alpha=(\mathcal{H}, y^k)$ is a reductive 
Whitehead move at 
$v = [\mathcal{H}, \underline{\underline{A}}]$ 
in 
$L(G)^F$ carried by $[\mathcal{H}, \underline{\underline{B}}]$
as described above, then there is some $a$ adjacent to
$k$ in the tree $T$ 
determined by $\underline{\underline{A}}$
such that $\alpha_a=(\mathcal{H}_a, x^k)$
is a reductive Whitehead move at 
$v_a = [\mathcal{H}_a, \underline{\underline{0}}]$ 
in $L(G_a)$ carried by 
$[\mathcal{H}_a, \underline{\underline{B_a}}]$.

Moreover, if $a_0, a_1, \ldots, a_m$ is a complete list of
vertices adjacent to $k$ in $T$ then
$$\alpha = (\alpha_{a_0})^{a_0} (\alpha_{a_1})^{a_1} \ldots 
(\alpha_{a_m})^{a_m}$$
and all of the (not necessarily reductive)
terms in the product commute.
\end{lem}

\begin{proof} The last assertion of the lemma
follows directly.  It remains to show that  at least one
$\alpha_{a_j}$ is reductive.  Let $a_0$ be the vertex
adjacent to $k$ in $T$ which is on the path from $k$ to $*$.
Let $a_1, \ldots, a_m$ be the other vertices adjacent to $T$,
ordered so that if $a_r < a_s$ in $X$ then $r < s$ as
integers.  Assume
as always that $y_k^k=1$.  Let $t$ be the least index such that
$y_j^k \not = 1$ for some $j \in I(a_t)$.  We show that
$\alpha_{a_t}$ is reductive.

Since $y_j^k \not = 1$ for some $j \in I(a_t)$, $\alpha$ must
change the norm of some letter in $G_{a_t}$.  By the minimality
of $t$, $\alpha$ does not change the norm of any letter
in $G_a$ for $a < a_t$.  Therefore, since $\alpha$ is 
reductive, it must be reductive on $G_{a_t}$.  
Lemma \ref{lem:frfgwordlength} now yields that $\alpha_{a_t}$
is reductive. 
\end{proof}

\begin{thm} \label{thm:frfgcontractible}
Let $F$ be a finite subgroup of $\Sigma Aut_1(G)$ that fixes
$[\mathcal{H}_0, \underline{\underline{A}}]$, suppose that
$[\mathcal{H}_0, \underline{\underline{A}}]$ is reduced in $L(G)^F$.
Then $L(G)^F$ is
contractible.
\end{thm}

\begin{proof}
We do this by induction, adding the ascending stars of 
nuclear vertices $v=[\mathcal{H}, \underline{\underline{A}}]$ in $L(G)^F$
step by step according to the norm 
of Proposition \ref{prop:frfgwellord}, always insuring that
the reductive part of the star $st(v)$ 
(of $v$ in $L(G)^F$) is contractible.

We follow the discussion of McCullough and Miller on pages 36-37
of \cite{mm}.  Namely, we first let $R_1$ be the reductive part
of the star of $v$.  We then let $R_2$ denote the full subcomplex of 
$R_1$ spanned by all vertices each of whose nontrivial based partitions
$\underline{\underline{B}}(i)$
is reductive.  
(Note that in the context of vertices in the star of 
$v=[\mathcal{H}, \underline{\underline{A}}]$ in $L(G)^F$,
a {\em trivial} based partition 
$\underline{\underline{B}}(i)$ is one which is 
equal to $\underline{\underline{A}}(i)$.)
Define a poset map $f_1: R_1 \to R_2$
as follows:  If $[\mathcal{H}, \underline{\underline{B}}]$ is
in $R_1$, then send it to 
$[\mathcal{H}, \underline{\underline{C}}]$, where each
nontrivial based partition
$\underline{\underline{B}}(j)$ with negative reductivity is
replaced by the trivial based partition with the same operative factor.
Since
$f_1(w) \leq w$ for all $w \in R_1$, Quillen's Poset Lemma \cite{q}
yields that $R_2$ is a deformation retract of $R_1$.

Let
$\underline{\underline{B}}(k)$ 
denote a partition corresponding
to some $[\mathcal{H}, \underline{\underline{B}}]$ in $R_2$
Suppose that $a_0, \ldots, a_m$ are the vertices adjacent to 
$k$ in the tree $T$ corresponding to $\underline{\underline{A}}$
(cf. the proof of Lemma \ref{lem:frfglargetosmall}.)
Then $\underline{\underline{B}}(k)$ is {\em admissible} if
each $(\underline{\underline{B}}_{a_j})^{a_j}(k)$ is either
trivial (that is, equal to $\underline{\underline{A}}(k)$)
or reductive.
Now let $R_3$ denote the full subcomplex of $R_2$ spanned by
all vertices each of whose nontrivial reductive based partitions 
is admissible.
Define a map $f_2: R_2 \to R_3$ by
combining all petals of
$\underline{\underline{B}}(k)$ containing elements of
$I(a_j)-\{k\}$
for each $j$ where  $(\underline{\underline{B}}_{a_j})^{a_j}(k)$ is 
nontrivial and not reductive.  As before, 
$f_2(w) \leq w$ for all $w \in R_2$ so that $R_2$
deformation retracts to $R_3$.

For each $a \in A$, let $R_2(st(v_a))$ denote the 
reductive portion of the star of $v_a$ in $L(G_a)$ where
where each based partition is either trivial
or reductive. Each nonempty $R_2(st(v_a))$ 
is contractible by Theorem \ref{thm:mmcontractible}.
Let $\bar A = \{a \in A: R_2(st(v_a)) \not = \emptyset\}$.
Recall that if $P_1$ and $P_2$ are posets, we can form their join
$P_1 \star P_2$ as the poset with elements 
$P_1 \cup (P_1 \times P_2) \cup P_2$. 
If $p_1, p_1' \in P_1$, $p_1 < p_1'$ in $P_1$
 $p_2, p_2' \in P_2$, and $p_2 < p_2'$ in $P_2$, then
in the poset $P_1 \star P_2$ we have
$(p_1, p_2) \geq (p_1', p_2')$, 
$(p_1, p_2) \geq p_1$,
$(p_1, p_2) \geq p_2$,
$p_1 \geq p_1'$, and
$p_2 \geq p_2'$.  This coincides with the more usual
definition of the join of two topological spaces
$$X \star Y = \frac{X \times [0,1] \times Y}
{(x,0,y) \sim (x,0,y'), (x',1,y) \sim (x,1,y)}$$
in the sense that the realization of
$P_1 \star P_2$ is homeomorphic to the join of the realization
of $P_1$ with that of $P_2$.
However, from Lemmas \ref{lem:frfgsmalltolarge}
and \ref{lem:frfglargetosmall} there is a poset isomorphism
$$f_3: R_3 \to \star_{a \in \bar A} R_2(st(v_a))$$
given by
$$f(\underline{\underline{B}}) =
\prod_{a \in \bar A, 
\underline{\underline{B_a}} \not = \underline{\underline{A_a}}}
\underline{\underline{B_a}}.$$
Since each poset in the join is contractible,
$\star_{a \in \bar A} R_2(st(v_a))$ is contractible.
\end{proof}


\end{document}